\documentclass[11pt]{article}
\usepackage{color}
\usepackage[toc,page]{appendix}
\usepackage[utf8]{inputenc}
\usepackage[margin=1in,left=1in]{geometry}
\usepackage{setspace}
\usepackage{graphicx}
\usepackage{amssymb, amsmath, amsthm, graphics}
\usepackage{mathabx,epsfig}
\usepackage[mathscr]{euscript}
\usepackage{tikz}
\usetikzlibrary{calc}
\usetikzlibrary{matrix}
\usepackage{hyperref}     

\usepackage{latexsym}
\usepackage[all]{xy}
\usepackage{color}

\usepackage{bbold}

\usepackage{tikz}
\input xy
\xyoption{all}
\usepackage{euscript}
\usepackage{multirow}
\usetikzlibrary{positioning}
\usetikzlibrary{shapes.geometric}
\usetikzlibrary{shapes.misc}
\usetikzlibrary{calc}
\usetikzlibrary{positioning}

\usepackage{pgflibraryarrows}
\usepackage{pgflibrarysnakes}

\usetikzlibrary{trees} 
\usetikzlibrary[trees] 

\newtheorem{theorem}{Theorem}
\newtheorem{definition}[theorem]{Definition}
\newtheorem{lemma}[theorem]{Lemma}
\newtheorem{corollary}[theorem]{Corollary}
\newtheorem{proposition}[theorem]{Proposition}

\newtheorem{remark}{Remark}

\newtheorem{example}{Example}

\numberwithin{equation}{section}

\renewcommand{\(}{\begin{equation*}}
\renewcommand{\)}{\end{equation*}}
\newcommand{\bea}{\begin{eqnarray*}}
\newcommand{\eea}{\end{eqnarray*}}

\def\proof {{Proof.}\hspace{7pt}}

\def\endofproof {\hfill{$\Box$}\\}

\newcommand{\beq}{\begin{equation}}
\newcommand{\eeq}{\end{equation}}

\newcommand{\tensor}{\otimes}

\newcommand{\op}[1]{\ensuremath{\operatorname{#1}}}

\numberwithin{equation}{section}

\renewcommand{\(}{\begin{equation}}
\renewcommand{\)}{\end{equation}}

\def\1{{\bf 1}}

\def\<{\langle}
\def\>{\rangle}

\def\Hom{{\rm Hom}}

\numberwithin{equation}{section}


\newcommand{\RR}{\ensuremath{\mathbb R}}

\newcommand{\ZZ}{\ensuremath{\mathbb Z}}
\newcommand{\Z}{\ensuremath{\mathbb Z}} 
\newcommand{\Q}{\ensuremath{\mathbb Q}}

\newcommand{\CC}{\ensuremath{\mathbb C}}

\author{Hisham Sati \footnote{hsati@nyu.edu}~
and Aliaksandra Yarosh\footnote{alexandra.yarosh@gmail.com}}
\title{Twisted Morava K-theory 
and connective covers of Lie groups}

\begin{document}
\maketitle

\begin{abstract}
Twisted Morava K-theory, along with computational techniques, including 
 a universal coefficient theorem and an Atiyah-Hirzebruch spectral sequence,
was introduced by Craig Westerland and the first author in \cite{SW}. 
  We employ these techniques to compute twisted Morava K-theory of all connective covers
   of the stable orthogonal group and stable unitary group,
  and    their classifying spaces, as well as spheres and Eilenberg-MacLane spaces.
    This extends to the twisted case some of the results of Ravenel and Wilson \cite{RW}
  and   Kitchloo, Laures, and Wilson \cite{KLW} for Morava K-theory. This also generalizes
  to all chromatic levels  computations by Khorami \cite{khorami} 
  (and in part those of Douglas \cite{Do})
  at chromatic level one, i.e. for the case of twisted K-theory.  
 We establish that for natural twists in all cases, there are only two possibilities: 
 either that the twisted  Morava homology vanishes, or that it is isomorphic to untwisted homology.
We also provide a variant on the twist of Morava K-theory, with mod 2 
cohomology in place of  integral cohomology.
\end{abstract}

\tableofcontents


\section{Introduction} 

Twisted cohomology and covering spaces have a very intimate relationship. 
  Twisted cohomology encodes additional data coming from a principal bundle on a space. 
 The earliest incarnation is perhaps cohomology with local coefficients which allows, in 
 particular, to define Poincar\'e duality for non-orientable manifolds, as long as we 
 consider the cohomology with coefficients given by the local system of the orientation 
 double cover. What we consider here is a generalization
  of the two concepts, namely that of a cover and that of a twisted cohomology, to generalized 
  covers and to generalized cohomology theories. For the former we will consider
  groups and their classifying spaces arising from cokilling homotopy groups 
  of the stable orthogonal and unitary groups, leading to their Whitehead towers. 
  On the latter, we will consider twisted Morava K-theory \cite{SW}, of 
  which twisted K-theory can essentially be viewed as a specific instance. 
    
    \medskip
    Morava K-theory $K(n)$ is in some sense an ``extension" of K-theory: it is a complex-oriented 
cohomology theory that defined for every integer $n$ and prime $p$, where $n$ is 
the chromatic level, i.e., the height of 
the corresponding formal group law. See 
\cite{JW}\cite{Wil}\cite{Wu}\cite{Ru}\cite{RW}\cite{GreenBook}\cite{HRW}

\noindent \cite{HKR}.   
Note that $K(n)_*(X)$ is always a coalgebra. When $X$ is an H-space, it 
is in addition a Hopf algebra. This will be the case for all the examples
for which we compute the Morava K-theory.
Therefore, their corresponding Morava 
K-homology theories will be Hopf algebras.

\medskip
The first author and Westerland \cite{SW}
 show that Morava K-theory $K(n)$, 
Morava E-theory $E_n$, and some of their variants admit twists by $K(\mathbb Z, n+1)$ bundles. 
The motivation for this theory came from string theory: it was conjectured by the first
author \cite{Sa09} that a twisted form of Morava K-theory and E-theory should describe 
an extension of the untwisted setting in his work with Kriz 
on anomalies at chromatic level two \cite{KS1}. A vast generalization 
of this conjecture is proved in \cite{SW}.

\medskip
The coefficients of the theory $K(n)_*=\ZZ/p[v_n , v_n^{-1}]$ form 
a graded field, which implies that $K(n)_*(-)$ always has a very useful computational 
tool manifested in the  K\"unneth isomorphism 
$$
K(n)_*(X \times Y) \cong K(n)_*(X) \otimes_{K(n)_*} K(n)_*(Y) .
$$
In our main examples, we do not quite have  products, but rather
fibrations. No similar general formula for calculations 
for fibrations  exits. However, when the fiber is an Eilenberg-MacLane 
space, a very useful tool is available: The universal coefficient theorem (UCT)
of Sati-Westerland \cite{SW}, which says that $K(\ZZ, n+2)$-twisted Morava
 K-theory $K(n)$ of a space is isomorphic to the untwisted Morava K-theory 
 of a certain $K(\ZZ, n+1)$-bundle over that space. This is a generalization of
a theorem by Khorami \cite{khorami} in the case of twisted K-theory, i.e.,
for $n=1$, where the fibers are $K(\Z, 2)\simeq \CC P^\infty$.

\medskip
The main examples that we consider arise from connective covers of classical Lie groups. 
The Whitehead tower is a way of approximating a space by a sequence of spaces of 
increasing connectivity that under mild assumptions converges to the original space. 
More precisely, it is a sequence of fibrations 
$X \langle m+1 \rangle \to X\langle m \rangle \to \dots \to X$, such that each 
$X \langle m \rangle$ is $(m-1)$-connected (i.e., the first non-trivial homotopy group 
occurs in degree $m$), and the induced map $\pi_i\left(X \langle m \rangle \right) \to \pi_i(X)$ 
is an isomorphism for $i \geq m$ (i.e., we are successively killing homotopy groups of $X$, 
starting from the bottom). In particular, $X\langle 2 \rangle$ is the universal cover of $X$. 
In general $X\langle m\rangle$ is called the \emph{$m$-connective cover} of $X$.
The main example appearing in our work is the connective covers of the stable orthogonal group 
and its classifying space. Some of the connective covers have distinguished names 
arising from string theory and M-theory \cite{SSS2}\cite{SSS3}\cite{9}

\vspace{-3mm}
\begin{equation}\label{eq:tower}
\xymatrix@R=1.2em{
&\vdots \ar[d]&\\
K(\ZZ, 11) \ar[r] &\op{BO}\langle 13\rangle = \op{BNinebrane} \ar[d]&\\
K(\ZZ/2, 9) \ar[r] &\op{BO}\langle 11\rangle = \op{B2\textit{-}Spin} 
\ar[d]\ar[r]^-{\frac{1}{240}p_3} &K(\ZZ, 12)\\
K(\ZZ/2, 8) \ar[r] &\op{BO}\langle 10\rangle = \op{B2\textit{-}Orient} \ar[d] \ar[r]^-{\alpha_{10}} 
&K(\ZZ/2, 10)\\
K(\ZZ, 7) \ar[r] & \op{BO}\langle 9\rangle = \op{BFivebrane} \ar[d] \ar[r]^-{\alpha_9} &K(\ZZ/2, 9)\\
K(\ZZ, 3) \ar[r] & \op{BO}\langle 8\rangle = \op{BString} \ar[d] \ar[r]^-{\frac{1}{6}p_2} &K(\ZZ, 8)\\
K(\ZZ/2,1) \ar[r] & \op{BO}\langle 4\rangle = \op{BSpin} \ar[d] \ar[r]^-{\frac{1}{2}p_1} &K(\ZZ,4)\\
K(\ZZ/2,0) \ar[r]  &\op{BO}\langle 2\rangle = \op{BSO} \ar[d]\ar[r]^-{w_2} &K(\ZZ/2,2)\\
& \op{BO} \ar[r]^-{w_1} &K(\ZZ/2,1)\;.
 }
\end{equation}
Notice that notation $\op{BO}\langle m \rangle$ means $(\op{BO})\langle m \rangle$ rather than 
$\op{B}(\op{O}\langle m \rangle)$.  In the above diagram, $w_1$ and $w_2$ are first and second 
Stiefel-Whitney classes, respectively, with the classes $\alpha_9$ and $\alpha_{10}$ their exotic 
generalizations  (see \cite{9}) and $\frac{1}{k}p_i$ is the fractional $i$th Pontrjagin class. 
Unlike $\op{SO}$ and $\op{Spin}$, which have classical descriptions as a Lie groups, $\op{String}$ 
-- as well as higher covers -- is not a Lie group.
There are various models for these groups in the literature, but we will not need such explicit 
characterizations in this paper, as we are mainly interested in the homotopy type.

\medskip
We compute twisted Morava K-theory of all the connective covers of $\op{BO}$ and $\op{BU}$
in Section \ref{Sec-BO}. 
We will consider various Eilenberg-MacLane (EM) fibrations arising in diagram \eqref{eq:tower} 
(and its complex counterpart), 
$$
\op{EM} \longrightarrow X\langle m\rangle \longrightarrow X ,
$$
to which we apply the UCT of \cite{SW},  leading to determination of the 
twisted Morava K-theory of $X$ from the untwisted Morava K-theory of 
$X\langle m \rangle$. The latter relies on recent extensive results of Kitchloo, 
Laures, and Wilson on Morava K-theory of spaces related to the classifying space 
$\op{BO}$ \cite{KLW}\cite{KLW2}, while the Morava K-theory of the fiber relies 
on the computations of Ravenel and Wilson of the Morava K-theory of Eilenberg-MacLane 
spaces \cite{RW}.

\medskip
The idea of calculating the twisted theory of a space via the untwisted theory of 
a related space also has other occurrences. Twisted K-theory of finite-dimensional 
projective spaces $\RR P^n$ was studied in \cite{BEJMS}, 
where they were shown to be computed, via T-duality, using the untwisted 
K-theory of the dual space (often a sphere). We will generalize the projective 
spaces to infinite dimensions and then to higher cohomological ranks,
i.e., Eilenberg-MacLane spaces.

\medskip
The instances where the fiber is a mod 2 Eilenberg-MacLane space require us to introduce twists by mod 2 
Eilenberg-MacLane spaces (of different rank), and consider a mod 2 analogue of the universal coefficient 
theorem from \cite{SW}, which we establish in Section \ref{Sec-mod2}. We observe that, 
in all cases considered, twisted Morava K-theory is either equal to untwisted Morava K-theory of these 
covers, or is zero altogether,  with a transition occurring after height 2. 


\medskip
The computations using the UCT of \cite{SW} rely on twisted homology rather than
on twisted cohomology versions of Morava K-theory. For $R$ a spectrum, 
it is generally easier to compute the twisted $R$-homology
rather than the twisted $R$-cohomology of spaces. 
Consider bundles over a space X with fiber the $R$-theory spectrum (satisfying certain conditions). 
Then twisted $R$-homology is naturally 
defined as the homotopy of the total space of the bundle, relative to the base
\cite{ABG} \cite{ABGHR}\cite{SW}.

\medskip
   Twisted K-theory  (see \cite{karoubi}\cite{Ros}\cite{BCMMS}\cite{AS1}\cite{AS2}) 
 twisted by $\op{PU}(\mathcal H)\simeq K(\ZZ, 2)$ bundles also has a homological counterpart,  
 namely twisted K-homology, discussed topologically in \cite{Do},
geometrically in \cite{Wa}, and analytically in \cite{Me}. 
Twisted geometric cycles for general CW-complexes is described in \cite{BCW}. 
Directly one prime at a time, a finite-dimensional (analytic) model 
for twisted mod $p$ K-theory by using twisted vectorial bundles with Clifford action  
is given in \cite{Go}. Morava K-theory at chromatic level one, $K(1)$, 
is one of $(p-1)$-summands in the Adams decomposition of K-theory at a prime $p$.
We will mainly be interested in the prime $p=2$, so that the relation between 
$K(1)$ and complex K-theory is even more direct.

  \medskip
 Twisted K-theory of compact Lie groups has been studied by physicists (see,
e.g., \cite{Br} \cite{MMS}\cite{GG}) as well as by mathematicians, 
starting with Douglas \cite{Do}.  The computation of the twisted K-groups
was extended to Lie groups which are not necessarily compact simple and simply 
connected in \cite{GG}\cite{MR}. The results for twisted K-theory $K^*(G, h)$,
for arbitrary choices of the twist $h$, are already rather complicated and hard to
understand. In \cite{MR} Mathai and Rosenberg introduced a new method for 
computing twisted K-theory using the Segal spectral sequence, giving simpler 
computations of certain twisted K-theory groups of compact Lie groups relevant 
for physics.  Twisted K-groups for Lie groups are shown to be trivial except in the 
simplest case of  $\op{SU}(2)$ and $\op{SO}(3)$.

\medskip
There has been a lot of activity in calculating the 
Morava K-theory of $BG$, for $G$ a finite group,  due to the Hopkins-Kuhn-Ravenel 
conjecture on evenness of Morava K-theory of $BG$ \cite{HKR} and its counterexamples by Kriz 
\cite{Kriz}. Note that all Eilenberg-MacLane spaces used in this 
paper have even Morava K-theory, by \cite{RWY}. Finite groups can support integral twists needed for 
Morava K-theory, as the integral cohomology of a $K(G, 1)$ 
space is in general nonzero, albeit all torsion. However, the 
homotopy degree of such a space is not supported by the computational 
techniques that we use here. It would be interesting to 
investigate twisted Morava K-theory of such spaces via alternative methods.

\medskip
Which chromatic level for Morava K-theory to use? The first effect we encounter here is {\it acyclicity}. 
From Ravenel-Wilson \cite{RW}, the $n$th Morava K-theory `sees' 
only the first $n$  Eilenberg-MacLane spaces. 
As far as localization goes, mod $p$ K-theory is the first in the 
Morava series $K(n)$. From Anderson and Hodgkin
\cite{AnH}, the space $K(G, n)$ is $K\Z/p$-acyclic
for all $n \geq 3$ and all $G$. This implies that, to the eyes of  
mod $p$ K-theory, the higher Eilenberg-MacLane spaces
might as well be points. 
For higher chromatic levels, Ravenel and Wilson \cite{RW} 
showed that for $K(n)$ at odd primes $p$, 
$$
\widetilde{K}(n)_{*} K(G, n+2)=0\;, \qquad 
\widetilde{K}(n)_{*} K(A, n+1)=0\;,
$$
for $A$ any finite abelian group, in particular a cyclic group $\Z/p$. 
Consequently, fibrations of connective covers stabilize in the sense that the 
$K(n)$-(co)homology would be the same for $X\langle m \rangle$ and 
$X \langle m+1\rangle$ after some critical $m$. 
Note that this is a general feature in the sense that every homology theory 
$E_*$ has a transitional dimension, at which all higher Eilenberg-MacLane 
spaces are $E_*$-acyclic (see \cite{Bo2}).

\medskip
For us this means that, 
in order to describe the connective covers involving EM spaces, we need
a chromatic level that is at least as high as the homotopy degree
of that space. That is, for String, BString, Fivebrane, and BFivebrane, 
we need at least $K(2)$, $K(3)$, $K(6)$ and $K(7)$, respectively. 
If we use lower chromatic levels then we will not be able to see the fiber in the 
fibration defining the connective covers. This then would mean that,
to the eyes of those Morava K-theories with corresponding chromatic
levels lower than the above threshold, the calculations essentially 
reduce to those of the total space, i.e., to Spin, BSpin, String and 
BString, respectively.

\medskip
The second effect that we witness is the {\it vanishing} of twisted
Morava K-theory of certain spaces. We consider Eilenberg-MacLane spaces in Section
\ref{Sec-EM}, as well as twists by fundamental classes of spheres in Section \ref{Sec-Sph}, 
and show that at height 2 and above the twisted Morava K-homology is  zero.
At first, it might seem surprising that the twisted homology would ever be zero at
all -- for untwisted homology, in the most trivial case of a point or contractible space one
has that the homology is equal to the coefficient ring. However, reduced \textit{twisted} 
homology of a point is zero, so vanishing twisted homology just means that the space 
behaves  like a point in our setting. For the case of twisted K-theory, it was shown
in \cite{MR} that the twisted K-groups of most Lie groups in fact vanish.

\medskip
Directly tracing the essence of the proofs of the theorems  in Section 
\ref{Sec-BO}, this method can be captured in the following general 
vanishing theorem for twisted Morava K-homology. 
  
  \newpage
\begin{theorem}
[Vanishing theorem for twisted Morava K-theory]
\label{thmvan} 
If a principal $K(\ZZ, n+1)$ bundle $\xi:E \to B$ is such that the induced map 
on Morava homology is a map of Hopf algebras, 
then composition with the Bockstein
map 
\[
K(n)_{*}\left(K(\ZZ/2, n) \right) 
\xrightarrow{\delta_*} K(n)_{*} (K(\ZZ, n+1)) \longrightarrow K(n)_{*} (E)
\]
is zero and hence twisted Morava homology of the base vanishes, $K(n)_*(B, \xi) = 0$.
\end{theorem}

Note that a map of Hopf algebra occurs, for example, 
 when $E \to B$ is a loop space map. 

\medskip
A third important phenomenon that we encounter is that of \emph{untwisting}. 
We show in Section \ref{Sec-BO} that for all connective covers, with the twist given by 
the corresponding class defining the fibration, the twisted Morava K-homology
is given by the underlying untwisted Morava K-homology.  There is, however, one notable exception
which, somewhat surprisingly, is the classifying space of a classical Lie group: we show in 
Proposition \ref{Prop20} that $K(2)$ of $B{\rm Spin}$ with a twist given by the generator 
$\tfrac{1}{2}p_1$ in fact vanishes.

\begin{theorem}
[Untwisting in the Whitehead tower]
\label{Thm2}
The twisted Morava K-homology of all groups in the Whitehead tower of 
the orthogonal and unitary groups, and their classifying spaces (except for $B{\rm Spin}$ 
in Proposition \ref{Prop20}),  with the canonical twist, is isomorphic to the underlying 
untwisted Morava K-homology. 
\end{theorem}

\medskip
This seems to be an instance of a more general phenomenon, and it does occur even in 
the case of K-homology. For $G$ a compact,  connected, simply connected, simple Lie 
group of rank $n$, the twisted K-homology ring with non-zero twisting class 
$k \in H^3(G; \Z)\cong \Z$ is an exterior algebra of rank 
$n-1$ tensored with a cyclic group \cite{Do}\cite{Br} 
$
K_*^{\tau (k)}(G) \cong \Lambda[x_1, \cdots, x_{n-1}] \otimes \Z/c(G, k)
$,
where $c$ is a constant which is, interestingly, most involved for the case of ${\rm Spin}(m)$.
This cyclic term can be viewed as the effect of the twist in the sense that one
gets the untwisted groups when $c$ is zero. This is certainly the case when the 
twist itself is zero, $k=0$. However, investigating the expressions for $c$ in 
\cite{Do} this can be zero for many nontrivial twists, hence effectively reducing the 
twisted theory to the untwisted one.  For the ``basic twist" $k=1$ one has that 
$c({\rm Spin}(4n-1), 1)$, $c({\rm Spin}(4n+1), 1)$, $c({\rm Spin}(4n+2), 1)$, and 
$c({\rm Spin}(4n), 1)$  are all equal to ${\rm gcd}\{1, 0\}$ and hence are zero.
This then implies that the cyclic factor is $\Z$, so that the twisted K-homology
reduces to the underlying exterior algebra. Therefore, for  the basic twist, one has 
$K_*^{\tau(1)}({\rm Spin}(n))\cong K_*({\rm Spin}(n))$. 
The same holds for instance for  other `special' values, e.g.  
$k=2n-1, 2n+1$ and $2n+3$. 
Douglas in \cite{Do} also finds that ${\rm Spin}(n)$ is very special in that, 
unlike other Lie groups, a detailed knowledge is needed of the twisted module 
structure on $\Z$ required to identify the cyclic orders.


\medskip
In addition to the universal coefficient theorem (UCT), another computational tool 
we will use is the twisted Atiyah-Hirzebruch spectral sequence (AHSS) from \cite{SW}, 
which approximates a twisted generalized (co)homology theory by usual (co)homology 
of successive quotients arising from nested filtrations of the underlying space $X$. This 
generalizes to higher chromatic levels the AHSS for twisted K-theory 
\cite{Ros1}\cite{Ros}\cite{AS1}\cite{AS2},
which in turn generalizes that of complex K-theory \cite{AH1}.

\medskip
Rationally, the integral version of Morava K-theory $K(n)$ with coefficients $K(n)_*=\ZZ[v_n , v_n^{-1}]$,
is isomorphic to $v_n$-periodic rational cohomology $K(n)^*(X) \otimes \Q \cong H^*(X; K(n)_*\otimes \Q)$, 
where $v_n$ is the generator of degree $2(p^n-1)$  (see \cite{KS1}\cite{SW}\cite{GS}).
The rational computations of those connective covers that we consider 
have been studied in \cite{SWh}. 
Here we consider one prime at a time, with the prime $p=2$ playing a distinguished 
role, hence concentrating instead on torsion information.

%
%

\medskip
This paper is organized as follows. We start in Section \ref{Sec-cl} by providing the 
setting and the main background that we need, 
recalling the results of Khorami on twisted K-homology in Section \ref{sec:khorami}, 
which we aim to generalize in later sections. 
In Section \ref{Sec-comp} we recall the fundamental computations in untwisted
Morava K-theory, mainly those of Ravenel and Wilson for Eilenberg-MacLane spaces
and Kitchloo, Laures and Wilson for spaces related to BO. Then in Section \ref{Sec-tmo}
we recall the main constructions of twists of Morava K-theory, as well as the 
computational tools developed in \cite{SW}, and which we will directly
apply in later sections. 
We also provide a partial characterization and computational tools for  
twisted $K(1)$, or mod 2 K-theory, in Section \ref{Sec-K1}. 
Then in the main parts of the paper in Section \ref{Sec-Comp} we present 
our computations and constructions. In Section \ref{Sec-BO} we compute the 
twisted Morava K-theory of those connective covers of the classifying spaces BO
and BU which have fibers integral Eilenberg-MacLane spaces. 
Then in Section \ref{Sec-EM} and Section \ref{Sec-Sph}, we provide the 
computations for twisted K-theory of  Eilenberg-MacLane spaces and of spheres, 
respectively. In order to complete the computation of connective covers from Section 
\ref{Sec-BO}, we compute the twisted Morava K-theory  of those covers 
in the tower \eqref{eq:tower} with 
fibers mod 2 Eilenberg-MacLane spaces in Section \ref{Sec-mod2}. This requires
a mild extension of some of the results of \cite{SW} to include mod 2 twists. 
Many of the computations in this paper have been part of the thesis of the 
second author \cite{Alex} under the guidance of the first author.

\vspace{.5cm}
\par{\large \bf Acknowledgements.} 
The authors are very grateful to Craig Westerland for tremendous help and encouragement 
throughout the course  of this project. His input certainly deserves that he be a co-author; 
however, he graciously declined, and the authors respect his decision. The authors also thank
Nitu Kitchloo for useful correspondence, and are grateful to the anonymous referee for the 
careful reading and excellent suggestions that have led to substantial improvements of the paper.

\section{Computations in Morava K-theory and (variations on) its twists} 
\label{Sec-cl}

\subsection{Computations in twisted K-homology}
\label{sec:khorami}

In this section we summarize the main results of Khorami \cite{khorami} which we will need. 
Among the results of this paper is a generalization of the main computations there, 
using the techniques from \cite{SW}. 
Consider a space $X$ with a twist $\tau \in H^3(X; \ZZ)$ of degree three
 and consider a $\mathbb{C} P^\infty$ principal bundle $\mathbb{C} P^\infty \to P_\tau \to X$ 
induced by $\tau$. 
This setup leads to a  spectral sequence (see \cite[Theorem 4.1]{EKMM})
$$
\op{Tor}_*^{K_*K(\ZZ, 2)}\big(K_*(P_\tau), K_*\big) 
\; \Longrightarrow \; \pi_* X^{T \tau}\cong K_*^\tau(X)\;.
$$
Even though $K_*$ is not a flat $K_*K(\ZZ, 2)$-module,  
the K-homology universal coefficient theorem says
\begin{theorem}[K-homology UCT \cite{khorami}]
\label{th:uct khorami}
 $$
K_*^\tau(X) \cong K_*(P_\tau) 
\otimes_{K_*(\mathbb{C}P^\infty)} \hat{K}_*\;,
$$
where $\hat{K}_*$ is just the coefficient ring $K_*$ with 
the $K_*(\mathbb{C}P^\infty)$-module structure obtained from the action 
of  $\mathbb{C} P^\infty$ on K-theory. 
\end{theorem}

The module structure is important for the purpose of this paper, and the following is a quick 
summary from how that works in the case of K-homology \cite{khorami}. The module structure 
for the more general case of twisted Morava K-theory in Section \ref{Sec-tmo} will be similar. 
The action of $\mathbb{C} P^\infty$  on $K$ is via tensor product with a complex line bundle 
${L}$, which is the main reason why twisted K-theory  is possible. The bundle $P_\tau$ admits 
a fiberwise action of  $\mathbb{C} P^\infty$ given by $\mathbb{C} P^\infty \times P_\tau \to P_\tau$.
The twisted K-homology of $X$ is given as \cite{khorami} 
$$
K_*^\tau (X) \cong K_*(P_\tau) \otimes_{K_*(\mathbb{C} P^\infty)} \hat{K}_* \;.
$$
The graded abelian group $\hat{K}_*$ is the same as $K_*$ with the $K_*(\mathbb{C} P^\infty)$-module 
structure being not the one obtained by collapsing $\mathbb{C} P^\infty$ to a point, 
$\mathbb{C} P^\infty \to {\rm pt}$ (cf. Remark \ref{rem-hop}),
but rather  coming from the map $\CC P^\infty \to {\rm GL}_1(K)$ and then using the 
tautological action of the target on $K_*$. 
For any principal $\mathbb{C} P^\infty$ bundle $P_\tau \to X$, K-homology
$K_*(P_\tau)$ is a $K_*(\mathbb{C} P^\infty)$-module, where the action of $\mathbb{C} P^\infty$ on the 
total space $P_\tau$ induces a map 
$
K_*(\mathbb{C} P^\infty \times P_\tau) \to K_*(P_\tau)
$.
Since $K_*(\mathbb{C} P^\infty)$ is free over the coefficients $K_*$, one has an
isomorphism 
$
K_*(\mathbb{C} P^\infty \times P_\tau) \cong K_*(\mathbb{C} P^\infty) \otimes_{K_*} K_*(P_\tau) 
$
which gives the module structure 
$K_*(\mathbb{C} P^\infty) \otimes_{K_*} K_*(P_\tau)  \to K_*(P_\tau)$.

\medskip
Note that the K-homology of $\mathbb{C} P^\infty$ can be given explicitly 
as follows (see
\cite{BlueBook}). From complex orientation, $K^*(\mathbb{C} P^\infty)= K^*({\rm pt})[[x]]$, 
where $x=L-1$, where $L$ is the Hopf line bundle over 
$\mathbb{C} P^\infty$. So there are unique elements 
$\beta_i \in K_{2i}(\mathbb{C} P^n)$, $1 \leq i \leq n$, such that 
$\langle x^k, \beta_i \rangle = \delta_i^k$, $1 \leq k \leq n$. 
The collection $\{\beta_0=1, \beta_1, \beta_2, \cdots \}$ forms a 
$K_*$-basis  for $K_*(\mathbb{C} P^\infty)$
$$
K_*(\mathbb{C} P^\infty)= K_*\{ \beta_0, \beta_1, \cdots\}=\ZZ[t, t^{-1}]\{\beta_0, \beta_1, \cdots\}\;.
$$
The main examples presented in \cite{khorami} are the following.
\begin{example} 
[Degree three integral Eilenberg-MacLane space $K(\mathbb{Z}, 3)$] 
 \label{ExKh1}
For a  twist $n: K(\ZZ, 3) \to K(\ZZ, 3)$, 
comparing the bundle $K(\ZZ, 2) \to P_n \to K(\ZZ, 3)$ 
with the path-loop fibration 
$K(\ZZ, 3) \to PK(\ZZ, 2) \simeq * \to K(\ZZ, 2)$
identities $P_n$ with $K(\ZZ/n\ZZ, 2)$. Then, invoking 
a result of Anderson and Hodgkin \cite{AnH} that 
$\widetilde{K}_*(K(\ZZ/n\ZZ, 2))=0$, gives the triviality 
of twisted K-homology $K_*^{(n)}(K(\ZZ, 3))=0$.
\end{example}

\begin{example}
[Three-Sphere $S^3$] 
    \label{ExKh2}
  Since $H^3(S^3; \ZZ)\cong \ZZ$ any twist  $n: S^3 \to K(\ZZ, 3)$
 corresponds to  an integer. 
The differential in the Atiyah-Hirzebruch-Serre 
spectral sequence is identified as
$d_3(\sigma_3)=n \beta_1$, where $\sigma_3$ is the natural
generator of $H_3(S^3, K_0(\CC P^\infty))$, corresponding to 
the natural generator of $H_3(S^3; \ZZ)$, and $\beta_1$ 
is the degree one generator of K-homology of $\CC P^\infty$ above, 
interpreted cohomologically as a map
$S^2=\CC P^1 \hookrightarrow \CC P^\infty$. 
The $\ZZ/2$-graded twisted K-homology is then 
\[
K_*^{(n)}(S^3)  \cong  (K_*(\CC P^\infty)/n\beta_1)
\otimes_{K_*(\CC P^\infty)} \hat{K}_*\cong K_*/n=\ZZ/n\ZZ\;,
\]
and vanishes for the basic twist $n=1$.
\end{example}

Example \ref{ExKh2} can also be deduced from the calculations of twisted K-homology of Lie groups 
\cite{Do}\cite{Br}\cite{MR} since $S^3 \cong \op{SU}(2)$. 
We will generalize the above two examples to higher dimensions 
and higher chromatic levels in Section \ref{Sec-EM} and Section \ref{Sec-Sph},
respectively. 

\subsection{Computations in Morava K-theory}
\label{Sec-comp}

We briefly describe some computational tools used in this paper.
First, just from the fact that $K(n)$ is complex-oriented, 
the Morava K-theory cohomology of classifying spaces of the unitary group
is given as \cite{BlueBook}\cite{GreenBook}

\vspace{-10mm}
\begin{eqnarray*}
    K(n)^*\left(\mathbb{CP}^{\infty}\right) &\cong& K(n)_{*}[[x]]\;, \\
    K(n)^*\left(\mathbb{CP}^{\infty}\times\mathbb{CP}^{\infty} \right) &\cong&K(n)_{*}[[x,y]]\;, \\
    K(n)^*\left(\op{BU}(n)\right) &\cong& 
    K(n)_{*}[c_1, c_2\dots c_n]\;,
\end{eqnarray*}    
    where $|x|=|y| =2$, and $|c_k| = 2k$.

\medskip
The fundamental computation here
 is the Morava K-theory of Eilenberg-MacLane spaces by Ravenel and Wilson. 
  While the Hopf algebra
  $K(n)^*K(\Z, n+1)$ is a power series algebra on one generator,
the dual Hopf algebra   
$K(n)_*K(\Z, n+1)$ turns out to be an algebra with generators $x$ satisfying 
$x^p=u x$ for a suitable unit $u$.  

\begin{theorem}[{\cite[Theorem~11.1, Theorem~12.1]{RW}}]
\label{th:RW}
Let $K(n)$ be Morava K-theory at prime $p$. 

\vspace{-2mm}
    \item {\bf (i)} $K(n)_*K(\ZZ/p^j,q) \cong K(n)_*$ for $q >n$. 
    \item {\bf (ii)} $K(n)_*K(\ZZ/p^j,n) \cong \bigotimes\limits_{i=0}^{j-1}R(a_i)$ \; 
    and  \;  $K(n)^{*}(K(\ZZ/p^j, n)) \cong K(n)_*[x]/x^{p^j} $,
    where the generator $x$ has dimension $|x| = 2\frac{p^n-1}{p-1}$, the element 
$a_k$ is dual to $(-1)^{k(n-1)}x^{p^k}$, and 

\vspace{-3mm}
$$
R(a_k) = \ZZ/p[a_k, v_n^{\pm 1}]/(a_k^p - (-1)^{n-1}v_n^{p^k}a_k)\;.
$$   
    
    \vspace{-.4cm}
    \item {\bf (iii)} $K(n)_*K(\ZZ,q+1) \cong K(n)_*$ for $q >n$. 
    \item {\bf (iv)} Let $\delta: K(\ZZ/2^j, q) \to K(\ZZ, q+1)$ be the Bockstein map, 
    and let $b_i := \delta_{*}(a_i)$. Then
    
    \vspace{-2mm}
    $$
    K(n)_*K(\ZZ, n+1) \cong \bigotimes\limits_{i=0}^{\infty}R(b_i) \qquad 
    \text{and} \qquad 
    K(n)^{*}(K(\ZZ/p^j, n)) \cong K(n)_*[[x]]\;.
    $$
\end{theorem}

\begin{remark}[Notation] 
    Note the difference in notation from \cite{RW}: Our $a_k$ and $b_k$ are originally $a_J$ 
    and $b_J$, with $J = (nk, 1, 2, \dots, n-1)$, and our $x$ is $x_S$ with $S= (1, 2, \dots , n-1)$.
 Note also that in \cite{RW} this theorem is only proven for odd primes, but it was later extended 
 to $p=2$ in \cite[~Appendix]{Wi}. See also \cite[Section 2.2]{HL} for more details.   
\end{remark}

The following example should help in understanding the 
relation between generators arising from $K(1)$, viewed
essentially as K-(co)homology, versus those arising from it 
being viewed as the first Morava  K-theory. 

\begin{example}
[$K(1)$ of $\CC P^\infty$] We set $p=2$. 
The formal group law arising from complex orientation gives that 
as an algebra $K(1)_{*} \CC P^\infty$ is generated by elements 
$
\beta_{(i)} \in K(1)_{2^i +1}(\CC P^\infty)
$,
with relations $\beta^{*2}_{(i)}=v_1^{2^i} \beta_{(i)}$, where $\beta_{(i)}:=\beta_{2^i}$ (see \cite{RW}). 
We keep in mind that in the comparison the lowest degree element
is $\beta_{(0)}:=\beta_1$. 

\vspace{-4mm}
{\small
$$
\begin{tabular}{|c||c|}
\hline
{\bf $K$-(co)-homology} & {\bf $K(1)$-(co)homology}
\\
\hline
\hline
$K^* (\CC P^\infty)=K^*({\rm pt})[[x]]=\Z[t^{\pm 1}][[x]]$, $|x|=2$
&
$K(1)^* (\CC P^\infty)=K(1)_*({\rm pt})[[x]]=
\Z/2[v_1^{\pm}][[x]]$, $|x|=2$
\\
\hline
$K_*(\CC P^\infty)=\Z[t^{\pm}]\{\beta_0, \beta_1, \cdots\}$, 
$\beta_j$  dual to $x^j$
&
$K(1)_*(\CC P^\infty)=\bigoplus_{k \geq 0} R(b_k)$, $|b_k|=2^{k+1}$,
$b_k$ dual to $x^{2^k}$
\\
\hline
\end{tabular}
$$
}
The indexing between the two columns in the table are related by
$j \leftrightarrow 2^k$. One consequence to keep in mind for bookkeeping
the indices of the generators  is that
$j$ starts from zero while $2^k$ starts from 1. 
In going from the first column to the second we have mod 2 reduction, from 
coefficients $\Z[t^{\pm 1}]$ to $\Z/2[v_1^{\pm 1}]$, where we also map 
$t \mapsto v_1$
in the process. 
%
%
\end{example}

When analyzing connective covers of $\op{BO}$ in Section \ref{Sec-BO}, 
we will need several results of Kitchloo, Laures and Wilson \cite{KLW}\cite{KLW2} 
about Morava K-theory of spaces related to $\op{BO}$. 

 \begin{theorem}[{\cite[Theorem~1.3]{KLW}}]\label{th:klw bo}
Let $\underline{\op{bo}}$, $\underline{\op{BO}}$, $\underline{\op{BSO}}$, 
$\underline{\op{BSpin}}$ denote the connective $\Omega$-spectra with zeroth spaces 
$\mathbb Z \times \op{BO}$, $\op{BO}$, $\op{BSO}$, and $\op{BSpin}$, respectively.
Let $E \to B$ be a connective cover with fiber $F$, and $B$ is one of the following: 
$\underline{\op{bo}}_i$, for $i \geq 4$, $\underline{\op{BO}}_i$, $\underline{\op{BSO}}_i$, 
$\underline{\op{BSpin}}_i$, for $i \geq 2$. Then the fibration $F \to E \to B$ induces the
following short exact sequence of Hopf algebras
\(
\label{eq:klw}
K(n)_{*} \longrightarrow K(n)_{*}(F) \longrightarrow K(n)_{*}(E) \longrightarrow
 K(n)_{*}(B) \longrightarrow K(n)_{*}\;, 
\)
where $K(n)$ is the Morava K-theory at prime $p=2$.
\end{theorem}

See Definition \ref{def-hopf} for short exact sequences of Hopf algebras. Several of 
the spaces used in the  the above theorem correspond to spaces in the Whitehead tower of the 
 orthogonal group \eqref{eq:tower}: for  $i=0$, we have
$$
\xymatrix@R=2em{
\underline{bo}_8=BO \langle 8 \rangle \ar[d] \ar[r] &
\underline{B{\rm Spin}}_0=B{\rm Spin} \ar[d] \ar[r] &
\underline{B{\rm SO}}_0=B{\rm SO}  \ar[d] \ar[r] &
\underline{B{\rm O}}_0=B{\rm O}\;,  \ar[d] &
\\
K(\Z, 8) & K(\Z, 4) & K(\Z/2, 2) & K(\Z/2, 1)
}
$$
and for $i=8$, we have 
$$
\xymatrix@R=2em{
\underline{B{\rm Spin}}_8=B{\rm Ninebrane} \ar[d] \ar[r] &
\underline{B{\rm SO}}_8=\op{B2\textit{-}Spin}  \ar[d] \ar[r] &
\underline{B{\rm O}}_8=\op{B2\textit{-}Orient} .   \ar[d] &
\\
K(\Z, 12) &  K(\Z/2, 10) & K(\Z/2, 9)
}
$$
See also Remark \ref{rem-BOO} below.

\begin{remark}
    It is worth noting that the maps in the short exact sequence  \eqref{eq:klw} are precisely the maps 
    induced by maps $F\to E\to B$ defining the connective cover. While it is not explicitly mentioned in the 
    statement of this theorem in \cite{KLW}, examination of the proof and the results upon which
    it relies  (\cite[Theorem ~4.2]{KLW}\cite[Proposition~2.0.1]{RWY}) makes it clear. This fact is 
    crucial to some of our computations in later sections.
\end{remark}

To deal with base spaces outside of the range specified by Theorem
\ref{th:klw bo}, we will need another exact sequence.

\begin{theorem}
[{\cite[Theorem~1.5]{KLW}}]\label{th:klw bspin}
    Let $K(n)$ be the Morava K-theory at $p=2$. There is a short exact sequence of Hopf algebras
    \(
\label{eq snake}
\xymatrix{
K(n)_{*} \ar@{->}[r] & 
K(n)_{*}K(\ZZ/2, 2) \ar@{->}[r]^-{\delta_*} & 
 K(n)_{*} K(\ZZ,3)  \ar@{->}[r] & 
 K(n)_* \op{BString} \ar@{->}`r[d]`[l]`[llld]_-{\delta}`[d][lld] & \\& 
 K(n)_* \op{BSpin}\ar@{->}[r] &
K(n)_* K(\ZZ, 4)\ar@{->}[r]^-{(\times 2)_*} &
K(n)_* K(\ZZ, 4) \ar@{->}[r] & 
K(n)_{*}
\;,}
\)
where $\delta_*$ is the map induced by Bockstein map $K(\ZZ/2, 1) {\longrightarrow} K(\ZZ,2)$, 
and $(\times 2)_*$ is the map induced by multiplication by $2$  on $K(\ZZ, 4)$. 
    \end{theorem}

A similar result holds for connective covers of the classifying space 
$\op{BU}$ of the unitary group $\op{U}$.

\begin{theorem}[{\cite[Section~2.6]{RWY},\cite[Theorem~1.2]{KLW}}]\label{th:bu}
Let $\underline{bu}$ denote the connective $\Omega$-spectrum with zeroth space $\op{BU}$, and 
let $E \to B$ is a connective cover with fiber $F$, and $B$ is  $\underline{bu}_i$, for $i \geq 2$.
Then the fibration $F \to B \to E$ induces the following short exact sequence of Hopf algebras:
    \[
    K(n)_{*} \longrightarrow K(n)_{*}(F) \longrightarrow K(n)_{*}(E) 
    \longrightarrow K(n)_{*}(B) \longrightarrow K(n)_{*}\;, 
    \]
    where $K(n)$ is Morava K-theory at a prime $p$.
\end{theorem}


\subsection{Twisted Morava K-theory}
\label{Sec-tmo}

 We will be particularly interested in twists of Morava K-theory 
 by Eilenberg-MacLane spaces. Given a map $X \to K(\Z, n+2)$ the induced bundle 
from the path-loop fibration over $K(\Z, n+2)$ in the diagram 
$$
\xymatrix@R=1.5em{
K(\Z, n+1) \ar[r] & P_H \ar[d] &&  PK(\Z, n+2)\simeq \ast  \ar[d] & K(\Z, n+1) \ar[l] 
\\
& X \ar[rr] && K(\Z, n+2) &
}
$$
is a principal $K(\Z, n+1)$-bundle, which is used in \cite{SW}  to twist $K(n)$-theory. 
 The twisting of $K(n)$ are constructed by studying the space of maps from 
 $K(\Z, m)$ to $\op{BGL}_1K(n)$. Using obstruction theory and showing that 
  the cobar spectral sequence collapses leads to the 
  space of twistings $tw_{K(n)}(K(\Z, n+2))$ being $\Hom_{K(n)_*-\op{alg}}
  (K(n)_*K(\Z, n+1), K(n)_*)$, hence is homotopically discrete. 
 A careful identification of this space leads to the following.

\begin{theorem}[{\cite[Theorem~3.1, 3.3]{SW}}]
\label{SW-tw}
   Consider twists of Morava K-theory $K(n)$ by Eilenberg-MacLane spaces $K(\mathbb Z, m)$. 
  
  \vspace{-2mm}
    \item {\bf (i)} There are no non-trivial twists for $m>n+2$.
     \vspace{-2mm}
    \item {\bf (ii)} For $m=n+2$, components of the space $Map(K(\mathbb Z, n+2), BGL_1 K(n))$ are contractible.
    \vspace{-2mm}
    \item {\bf (iii)} If $p \neq 2$ then the set of twists $\op{tw}_{K(n)}(K(\mathbb Z, n+2)) \cong \op{tw}_{K(n)}(*)$.
     \vspace{-2mm}
    \item {\bf (iv)} If $p = 2$, then the set of twists is a group isomorphic to dyadic integers, i.e., one has \newline
    $\op{tw}_{K(n)}(K(\mathbb Z, n+2)) \cong \mathbb Z_2$. 
  \end{theorem}
Because of this theorem,  henceforth we restrict our attention to the case $p=2$ when dealing with Morava K-theory
twisted by integral Eilenberg-MacLane spaces.
We will also need the following definition from \cite{SW}, of which we will provide a mild 
generalization in Section \ref{Sec-mod2}.

\begin{definition}
\label{def-u}
{\bf (i)} 
    A \emph{universal twist} $u$ is an element of $\op{tw}_{K(n)}(K(\mathbb Z, n+2)) \cong \ZZ_2$ 
    that is a topological generator. 
         \vspace{-2mm}
\item {\bf (ii)}
    Given a class $H \in H^{n+2}(X;\ZZ) \cong [X, K(\ZZ, n+2)]$, the 
    \emph{$H$-twisted Morava K-theory} is defined as
    \[K(n)_{*}(X; H) := K(n)^{u(H)}_{*}(X), \]
    and similarly for cohomology.
\end{definition}

Once and for all we make an arbitrary fixed choice of $u$.
Defined this way, twisted Morava K-theory has all the properties that we would like it to have
{\cite[Theorem~4.1]{SW}}, including 
normalization, namely that if $H=0$ then $K(n)^*(X; H)=K (n)^*(X)$,
and the module property, i.e., $K(n)^*(X; H)$ is a module over $K(n)^0(X)$.
\footnote{Note that here the notation with the superscript $0$ means the untwisted Morava K-theory of $X$,
 \emph{not} its degree 0 part.}

\medskip
Notice that any cohomology class in $H^{n+2}(X; \ZZ)$ can be interpreted as a $K(\ZZ, n+1)$-bundle. The main 
computational tool we employ is the relationship between twisted homology of the base and untwisted
cohomology of the total space.

\begin{theorem}[Universal coefficient theorem {\cite[Theorem~4.3]{SW}}]\label{th:uct morava}
       Let $H \in H^{n+2}(X)$, and $P_H \to X$ be the principal $K(\mathbb Z, n+1)$ bundle over $X$,
       classified by $H$. Then
 \[
K(n)_*(X;H) \cong K(n)_*(P_H) \otimes_{K(n)_*(K(\mathbb Z, n+1))} K(n)_*\;.
 \]
    Here $K(n)_*(P_H)$ is a $K(n)_*(K(\mathbb Z, n+1))$ module since $P_H$ is a $K(\ZZ,n+1)$ bundle,
    and $K(n)_*$ is made into $K(n)_*(K(\mathbb Z, n+1))$ by sending $b_0$ to $v_n$ and $b_i$ to 0 for all $i>0$, 
 where we make use of Theorem \ref{th:RW} for the structure of $K(n)_*(K(\mathbb Z, n+1))$.
\end{theorem}

The module structure here is also crucial, as in the
 case of chromatic level one, i.e. twisted K-homology, 
 emphasized after Theorem \ref{th:uct khorami}. 
 Since this is completely analogous with obvious changes,
 we will not repeat the discussion. 
 
 \medskip 
 Since $K(n)$ of Eilenberg-MacLane spaces is known by the work of 
 Ravenel-Wilson (Theorem \ref{th:RW} above), 
Theorem  \ref{th:uct morava} theoretically reduces the problem of computing twisted homology to computing homology of the total space.

\medskip
Another very useful computational 
tool at our disposal is the Atiyah-Hirzebruch spectral sequence (AHSS)
for twisted Morava K-theory. 
 
\begin{theorem}[{\cite[Theorem~5.1]{SW}}]\label{th:ahss}
        For $H \in H^{n+2}(X)$, there is a spectral sequence converging to twisted Morava K-theory 
    $$
    E_2^{p, q} = H^p(X, K(n)^q) \Longrightarrow K(n)^*(X; H)\;.
    $$
    The first possible nontrivial differential is $d_{2^{n+1}-1}$, given by
    \[
    d_{2^{n+1}-1}(x) = \big(Q_n(x) + (-1)^{|x|}  x\cup (Q_{n-1} \cdots Q_1(H))\big)\;.
    \]
\end{theorem}

Here $Q_n$ is the cohomology operation known as $n$th Milnor primitive at the prime $2$. It may be 
defined inductively  as $Q_0 = Sq^1$, and  $Q_{j+1}=Sq^{2^j}Q_j - Q_j Sq^{2^j}$, where 
$Sq^j: H^n(X;\ZZ/2) \to H^{n+j}(X;\ZZ/2)$ is the $j$-th Steenrod square
operation in mod 2 cohomology.

\subsection{Twisted $K(1)$ and twisted mod $2$ K-theory}
\label{Sec-K1}

The `default' coefficients for a (co)homology theory is the integers. 
We are also interested in other coefficients, mainly the integers mod $p$.
Morava K-theory $K(n)$ is defined one prime at a time, while K-theory $K$
does not depend on a prime. In order to be able to compare the latter to the former 
at chromatic level $n=1$, we need to restrict K-theory to seeing 
one prime at a time. 

\medskip
Given a spectrum $E$ and any abelian group $G$, one can define a 
\emph{spectrum with coefficients in $G$} (and, correspondingly, $E$-cohomology 
with coefficients in $G$) as the smash product
$
EG : = E \wedge SG
$,
where $SG$ is a Moore spectrum of type $G$. This spectrum $EG$ satisfies the 
following short exact sequences \cite[~p.200]{BlueBook}
\[
\xymatrix{
0\ar[r] & \pi_n(E)\tensor G \ar[r] & \pi_n(EG) \ar[r] & \op{Tor}_1(\pi_{n-1}(E), G) \ar[r] & 0\;, 
}
\label{eq:coefficients ses homotopy}
\]
or more generally, for any space $X$, the following holds at the level of 
$E$-homology
\(
\label{seq-EG}
\xymatrix{
0\ar[r] & E_n(X)\tensor G \ar[r] &  (EG)_n(X) \ar[r] & 
\op{Tor}_1(E_{n-1}, G) \ar[r] & 0\;,
}
\)
where $(EG)_{n}$ 
 denotes the homology theory
 corresponding to the spectrum $EG$. 

  \medskip
        Taking $G = \ZZ/p$, we get ``mod $p$" $E$-homology. 
%
%
%
For $E=K$, the complex K-theory spectrum, we can define K-theory 
with mod $p$ coefficients as a spectrum $K\ZZ/p = K \wedge S\ZZ/p$, 
where $S\ZZ/p$ is the Moore spectrum of type $\ZZ/p$. Moreover, there is 
the short exact sequence \eqref{seq-EG}, in this case 
relating $K\ZZ/p$ homology and $K$-homology of a space $X$.
On the other hand, by a classical result of Adams \cite{BlueBook}, mod $p$ K-theory decomposes into a 
sum of $p-1$ successive suspensions of $K(1)$. In our case of interest, $p=2$, mod 2 K-theory 
coincides with $K(1)$ since there is only one summand, and so $K(1)$ fits into the exact sequence 
\[
0 \longrightarrow  K_n(X)\otimes \ZZ/2 \longrightarrow 
 K(1)_n(X) \longrightarrow  \op{Tor}_1\big(K_{n-1}(X), \ZZ/2\big) \longrightarrow 0\;.
 \]
We would like to establish a twisted version of this relationship.

\begin{theorem}
[Exact sequence for twisted mod $p$ K-theory]
\label{th:mod p}
Let $X$ be a topological space and $\tau_3 \in H^{3}(X;\ZZ)$. 
Then we have the following exact sequence 
   \(
   \label{ES-KnHZ2}
    0 \longrightarrow  K_n^{\tau_3}(X) \otimes \ZZ/2
    \longrightarrow  K(1)_n(X;\tau_3) \longrightarrow 
     \op{Tor}_{1}\big(K_{n-1}^{\tau_3}(X), \ZZ/2\big) \longrightarrow 0\;, 
    \)
where $K^{\tau_3}_*(X)$ is K-theory twisted by $\tau_3=\tau$, as defined in 
Section \ref{sec:khorami}.
\end{theorem}
\proof
 Let $P$ be the total space of the principal $K(\ZZ,2)$ bundle classified by a degree three class
 $\tau_3 \in H^{3}(X;\ZZ) \cong [X, K(\ZZ,3)]$. The Thom spectrum corresponding to the twisted theory  
 $K(1)(X;H)$ is 
 $$
 X^{u(H)} \simeq \Sigma^{\infty}P_+ \wedge_{\Sigma^{\infty}K(\ZZ, 2)} K(1) \simeq
 \Sigma^{\infty}P_+ \wedge_{\Sigma^{\infty}K(\ZZ, 2)} (K \wedge S\ZZ/2)\;,
 $$ 
 where $S\ZZ/2$ is the Moore spectrum of type $\ZZ/2$, as above.
     Now $K(\ZZ, 2)$ acts trivially on $S\Z/2$ and the associativity of a ``mixed" smash product is established via \cite[Proposition~3.4]{EKMM}, 
   so that $X^{u(\tau_3)} \simeq \left(\Sigma^{\infty}P_+ \wedge_{K(\ZZ,2)} K\right)\wedge S\ZZ/2$. 
   On the other hand, from \cite[Section 7.2]{ABG}, the Thom spectrum for K-theory twisted by $\tau_3$ is 
    exactly $\Sigma^{\infty}P_+ \wedge_{K(\ZZ,2)} K \simeq X^{T(\tau_3)}$, where $T: K(\ZZ,3) \to K$ 
    is the map defined in Section \ref{sec:khorami}. Hence $X^{u(\tau_3)} \simeq X^{T(\tau_3)} \wedge S\ZZ/2$.
    Then, using the exact sequence \eqref{seq-EG}, we obtain
    \[
    0 \longrightarrow \pi_n\big( X^{T(\tau_3)}\big) \otimes \ZZ/2
    \longrightarrow \pi_n\big(X^{u(\tau_3)}\big) 
    \longrightarrow\op{Tor}_{1}\big(\pi_{n-1}\big( X^{T(\tau_3)}\big), \ZZ/2\big) \longrightarrow 0\;. 
    \]
 But we defined twisted homology precisely as homotopy groups of the Thom spectrum. 
 Therefore, we can rewrite this exact sequence as \eqref{ES-KnHZ2}.
\endofproof

We will make use of this theorem in Section \ref{Sec-Sph} when calculating 
the twisted Morava K-theory of spheres leading to Proposition \ref{thm 3-sph}. 
We will also consider a variation in Section \ref{Sec-mod2}, where the twists
themselves will take values in cohomology with $\ZZ/2$ coefficients.


%

\section{Computations of twisted Morava K-homology of spaces}
\label{Sec-Comp}

\subsection{Twisted Morava homology of connective covers of $\op{BO}$ and $\op{BU}$}
\label{Sec-BO}
In this section we will apply our computations to the main 
examples, which are connective covers of the stable orthogonal and 
unitary groups (see diagram \eqref{eq:tower}).
Recall that we are working with Morava K-theory at $p=2$
for the rest of the paper (except for the last section).

\medskip
From \cite{RW} (Theorem \ref{th:RW} above) we have that $K(n)_*(K(\Z, j))$ is trivial for 
$j> n+1$. With $n=1$, it follows that $K(1)_*(K(\Z, j>2)$ is trivial. 
Therefore, we  immediately 
 get that all connective covers of Spin, that is, String, Fivebrane etc. have the same
twisted $K(1)$-homology, while all classifying spaces covering $\op{BSpin}$, that is, 
$\op{BString}$,
$\op{BFivebrane}$, etc., have the same twisted $K(1)$-homology.
More precisely, we have   the following. 

\begin{lemma} 
[Twisted $K(1)$-homology of connective covers] 
\label{LemmaK1}
For $m> 3$, let ${\rm Spin}\langle m \rangle$ 
be a connective cover of Spin and ${\rm BSpin}\langle m \rangle$ be the 
 connective covers of  the corresponding classifying space BSpin. Then, for any $m'>m$, 
 we have that
\vspace{-2mm}
\item {\bf (i)} The map ${\rm Spin}\langle m' \rangle \to \op{Spin} \langle m \rangle$
 is an isomorphism in twisted $K(1)_*$. 
 
 \vspace{-2mm}
\item {\bf (ii)} The map ${\rm BSpin}\langle m' \rangle \to \op{BSpin}\langle m \rangle$
 is an isomorphism in twisted $K(1)_*$. 
 
  In fact, in both cases the homology is necessarily 
 untwisted.
 \end{lemma}

\proof
From \cite[Lemma 3.3]{HRW} one has that fibrations $F \to E \to B$ with $F$ $K(n)$-acyclic induce
an isomorphism $K(n)_*(E) \overset{\simeq}{\to} K(n)_*(B)$. By inspection, from the Whitehead tower 
\eqref{eq:tower}, the fibers of Fivebrane and beyond begin with $K(\Z, 6)$ through
$K(\Z/2, 7)$, $K(\Z/2, 8)$, and $K(\Z, 10)$ etc., all of which are $K(1)$-acyclic.
Similarly for the fibers of the corresponding classifying spaces.  This proves the lemma in the 
untwisted case. 
The twisted case follows trivially by noticing that,  for $m>3$, $H^3({\rm Spin}\langle m \rangle; \Z)$
and $H^3(\op{BSpin}\langle m \rangle ; \Z)$ are both 0,
yielding trivial twists.  
%
%
\endofproof

The omitted case in Lemma \ref{LemmaK1}(i), namely $m=3$,
corresponding to String, admits nontrivial twists. It is 
delicate and in the unstable case would use a careful  analysis 
of twisted K-homology of the Lie group ${\rm Spin}(q)$ according to rank, as in 
\cite{Do}. 
With the mod 2 version of Khorami's Theorem \ref{th:uct khorami} 
(as in Section \ref{Sec-mod2}), we  have 
$$
K(1)_*^\tau ({\rm Spin}(q)) \cong K(1)_*({\rm String}(q)) \otimes_{K(1)_*K(\Z, 2)} K(1)_*\;.
$$
Through this one could use the results in \cite{Do} 
(also \cite{Br}\cite{MMS}\cite{GG}\cite{MR})
to potentially gain some understanding 
of $K(1)_*({\rm String}(q))$. However, we will not do so here, 
although the left hand side has been calculated explicitly for the specific case of 
$q=3$ in Proposition \ref{thm 3-sph} below.

\begin{remark}[Equivalences in higher degrees] 
We have a similar pattern for  higher chromatic levels which are less than   
the degree of the Eilenberg-MacLane
fiber, i.e., 
$B{\rm String}$ for $n<2$, Fivebrane for $n<6$, and $B{\rm Fivebrane}$ for
$n<7$. Therefore, in order to get something nontrivial we have to 
consider $K(n\geq 3)(B{\rm String})$, $K(n\geq 7)({\rm Fivebrane})$,
and $K(n\geq 8)(B{\rm Fivebrane})$.
\end{remark}

\medskip
Note that this is somewhat related to \cite{SWh}. There it was shown that,
rationally, all higher connective covers can be described using Spin rather
than having to define a structure at level $r$ from that at level $r-1$. 
What we have is an analogue of rationalization, in the sense that Morava
$K(n)$-homology only sees Eilenberg-MacLane spaces of certain degrees,
and the ones not seen can be viewed as rationalized as far as $K(n)$-homology
goes. Indeed, for a fibration $F \to E \overset{j}{\to} B$  with 
fiber $F$ which is $K(n)_*$-acyclic, the map $j$ is a $K(n)_*$-equivalence (see 
\cite[Lemma 3.3]{HRW}). From \cite{RW}, $K(\Z, i)$ is $K(n)_*$-acyclic
if and only of $i> n+1$. Consequently, we immediately have the following. 

\begin{lemma}
[$K(n)$-equivalences for connective covers] We have the following

\vspace{-2mm} 
\item {\bf (i)} $K(\Z, 3)$ is $K(1)_*$-acyclic so that $B{\rm Spin}$ and $B{\rm String}$ are
$K(1)_*$-equivalent.

\vspace{-2mm} 
\item {\bf (ii)} $K(\Z, 6)$ is $K(n)_*$-acyclic for $n<5$ 
 so that String and Fivebrane are
$K(n<5)_*$-equivalent.

\vspace{-2mm} 
\item {\bf (iii)} $K(\Z, 7)$ is $K(n)_*$-acyclic for $n<6$ 
so that $B{\rm String}$ and $B{\rm Fivebrane}$ are
$K(n<6)_*$-equivalent.
\end{lemma}

Similarly one can go up the Whitehead tower \eqref{eq:tower} 
and deduce analogous statements  for 
2-Orient, 2Spin and Ninebrane, etc. This is compatible with Bousfield's result 
that for any space
$X$, each $K(n)_*$-equivalence of spaces if a $K(m)_*$-equivalence
for $1 \leq m \leq n$ \cite{Bo-Mo}.

\medskip
Our main result in this section is the following statement in the scope
 of Theorem \ref{Thm2}.

\begin{theorem}
[Twisted Morava K-homology of $\op{BO}\langle k \rangle$]
\label{Kn-BO}
Let $K(\mathbb{Z}, n+1) \to \op{BO}\langle n+3 \rangle \to \op{BO}\langle n+2 \rangle $ be 
a fibration defining a connective cover of $\op{BO}$ (so $n = 2 \mod\, 8 $ or $n= 6 \mod\, 8$) 
and $n \geq 6 $, and let $\tau_{n+2}$ denote the corresponding class in 
$H^{n+2}\left(\op{BO}\langle n+2\rangle; \mathbb Z \right)$. Then the twisted Morava K-homology 
of the classifying space $\op{BO}\langle n \rangle$ and the corresponding loop group 
$\op{O} \langle n-1 \rangle:=\Omega \op{B O} \langle n \rangle$ are given, respectively, as 
\begin{eqnarray*}
K(n)_*(\op{BO}\langle n+2 \rangle; \tau_{n+2}) &\cong& K(n)_*(\op{BO}\langle n+2 \rangle)\;,
\\
K(n-1)_*(\op{O}\langle n+1 \rangle; \tau_{n+1}) &\cong& K(n-1)_*(\op{O} \langle n+1 \rangle)\;,
\end{eqnarray*}
where $\tau_{n+2}$ is the twisting class and $\tau_{n+1}$ is its looping. 
\end{theorem}
 Our main tool in this section, which goes towards proving the above theorem,
is the exact sequence of Kitchloo-Laures-Wilson, Theorem \ref{th:klw bo} above,
where $E \to B$ is a connective cover with fiber $F$, and $B$ is one of the following: 
$\underline{\op{bo}}_i$, for $i \geq 4$, $\underline{\op{BO}}_i$, 
$\underline{\op{BSO}}_i$, $\underline{\op{BSpin}}_i$, for $i \geq 2$.
  It is worth noting that the maps in this short exact sequence are precisely the maps 
  induced by maps $F\to E\to B$ defining the connective cover. 
%
This fact is crucial to our computations.

\begin{remark}
\label{rem-BOO}
Notice that if $\underline{\op{bo}}$, $\underline{\op{BO}}$, $\underline{\op{BSO}}$, 
$\underline{\op{BSpin}}$ denote the connective $\Omega$-spectra with zeroth spaces $\mathbb Z \times \op{BO}$, $\op{BO}$, $\op{BSO}$, and $\op{BSpin}$, then we have the following equivalences
$$
\op{BO}\langle 8k \rangle \simeq \underline{\op{bo}}_{8k}, \qquad
\op{BO}\langle 8k+1 \rangle \simeq \underline{\op{BO}}_{8k}, \qquad
\op{BO}\langle 8k+2 \rangle \simeq \underline{\op{BSO}}_{8k}, \qquad
\op{BO}\langle 8k+4 \rangle \simeq \underline{\op{BSpin}}_{8k}.
$$
In particular, $\op{BString} = \op{BO}\langle 8 \rangle = \underline{\op{bo}}_8$, $\op{BFivebrane} = \op{BO}\langle 9 \rangle = \underline{\op{BO}}_8$ and, since the spectra in question are $\Omega$-spectra, $\op{String} = \Omega \op{BString} \simeq \Omega \underline{\op{bo}}_{8} \simeq \underline{\op{bo}}_7$ and 
$\op{Fivebrane} \simeq \underline{\op{BO}}_7$.
\end{remark}

%

We will also need to use some basic facts about Hopf algebras. Standard references include 
\cite{hopf}\cite{MM}\cite{Sw}. 
From an algebraic point of view, exact sequences and 
extensions of Hopf algebras are described, for instance, in \cite{AD}.

\begin{definition}
\label{def-hopf}
Let $A,B,C$ be commutative Hopf algebras over a field $k$. Suppose $i: A \to B$ is an injection of Hopf algebras, and
$j: B \to C$ is a surjection of Hopf algebras. Then if $j$ induces the isomorphism 
$C \cong B/i(A^{+})B$ of Hopf algebras, where $A^{+}$ denotes 
the augmentation ideal of $A$, we say that 
\[
k \to A \xrightarrow{i} B \xrightarrow{j} C \to k
\] 
is a \emph{short exact sequence of Hopf algebras}.
\end{definition}

\begin{remark}[Hopf algebras and module structure]
\label{rem-hop}
If we have an injective Hopf algebra morphism $i:A \to B$, we can view $B$ as an $A$-module, and then 
$B/i(A^{+})B \cong B \tensor_A k$ (sometimes also denoted $B/\!/A$). 
Therefore, in \textit{any} short exact sequence of commutative $k$-Hopf algebras 
$k \to A \to B \to C \to k$, we always have $C \cong B \tensor_A k$.
In particular, in the exact sequence \eqref{eq:klw}, we have
$K(n)_*(B) \cong K(n)_*(E) \otimes_{K(n)_*(F)} K(n)_*$, where the 
$K(n)_*(F)$-module structure on $K(n)_*$ is in this case via the mpa $F \to {\rm pt}$ 
(cf. the discussion after Theorem \ref{th:uct khorami}). 
\end{remark}

Recall also the following classical result.

\begin{theorem}[{\cite[~Th. 4.4]{MM}}]\label{MMt1}
 If $A$ is a connected Hopf algebra over a commutative ring with unity $K$, $B$ is a connected $A$-module 
 coalgebra, $i: A \to B$, $\pi: B \to K\tensor_A B$ are the canonical morphisms, and the sequences 
 $0 \to A \xrightarrow{i} B$, $B \xrightarrow{\pi} K\tensor_A B \to 0$ are split exact as sequences 
 of graded $K$-modules, then there exists $h: B \to A \tensor_K (K\tensor_A B)$, which is 
 an isomorphism of $A$-modules.
\end{theorem}

If we take the ring $K$ to be $K(n)_*$, then any module over $K(n)_*$ is free, and so any exact 
sequence is split automatically  as $K(n)_*$-modules.

\medskip
\emph{Proof of Theorem \ref{Kn-BO}}.
Now consider the exact sequence \eqref{eq:klw} again. Applying Theorem \ref{MMt1} with $A = K(n)_*(F)$ 
and
\footnote{Note the multiple uses of the notation $B$: as the base of the fibration as
well as the middle Hopf algebra. We hope the distinction will be clear from the context.}  
 $B = K(n)_*(E)$, we obtain an isomorphism
\[
K(n)_*(E) \cong K(n)_*(F) \otimes_{K(n)_*} \left(K(n)_* \otimes_{K(n)_*(F)} K(n)_*(E) \right).
\]
But the latter term is precisely $K(n)_*(B)$, as mentioned above, so that we have
the isomorphism
$
K(n)_*(E) \cong K(n)_*(F) \otimes_{K(n)_*} K(n)_*(B)
$.
Now, if $M$ is any $K(n)_*(F)$-module, we have
\begin{align*}
M \otimes_{K(n)_*(F)}K(n)_{*}(E)  &\cong M \otimes_{K(n)_*(F)} K(n)_*(F) \otimes_{K(n)_*} K(n)_*(B)\\
&\cong M\otimes_{K(n)_*}K(n)_*(B),
\end{align*}
as $K(n)_*$-modules. We emphasize that here the actual $K(n)_*(F)$-module 
structure of $M$ is irrelevant. 

Now take $F \to E \to B$ to be the connective cover $K(\mathbb{Z}, n+1) \to 
\op{BO}\langle n+3 \rangle \to \op{BO}\langle n+2 \rangle $ for $n = 2$ mod $8$ 
or $n = 6$ mod $8$. If $n = 2$ mod $8$, then $n+2 = 4$ mod $8$, and $\op{BO}\langle n+2 \rangle \cong \underline{\op{BSpin}}_{n-2}$, by Remark \ref{rem-BOO}. If  $n = 6$ mod 
$8$, then $n+2 = 0$ mod $8$ and, by the same Remark, 
$\op{BO}\langle n+2 \rangle \cong \underline{\op{bo}}_{n+2}$. So, as long as 
$n \geq 6$, $\op{BO}\langle n+2 \rangle$ is one of the spaces that can serve as the base 
for the fibration in the sequence \eqref{eq:klw}. So, in particular, for any 
$K(n)_*(K(\ZZ, n+1))$-module $M$,
\begin{equation}
\label{eq:m}M \otimes_{K(n)_*K(\ZZ, n+1)} K(n)_{*}(\op{BO}\langle n+3\rangle) 
\cong M\otimes_{K(n)_*}K(n)_*(\op{BO}\langle n+2\rangle)\;. 
\end{equation}
The UCT (Theorem \ref{th:uct morava}) for the bundle 
$K(\mathbb{Z}, n+1) \to \op{BO}\langle n+3 \rangle \to \op{BO}\langle n+2 \rangle $ 
states that
\[
K(n)_*(\op{BO}\langle n+2 \rangle; \tau_{n+2}) \cong
K(n)_*(\op{BO}\langle n+3 \rangle) \otimes_{K(n)_*K(\ZZ, n+1)} K(n)_* \;,
\]
with a special $K(n)_*K(\ZZ, n+1)$-module structure on the latter factor from 
Theorem \ref{th:uct morava}. But taking 
$M = K(n)_*$ in \eqref{eq:m}, we see that

\vspace{-9mm}
\begin{align*}
    K(n)_*(\op{BO}\langle n+2 \rangle; \tau_{n+2}) &
    \cong K(n)_*\otimes_{K(n)_*} K(n)_*(\op{BO}\langle n+2 \rangle) 
    \\
    &\cong K(n)_*(\op{BO}\langle n+2 \rangle).
    \end{align*}
    
    \vspace{-1cm}
\endofproof

We highlight that Theorem \ref{Kn-BO} indicates that for the natural twists 
associated with connective covers of the orthogonal group and their classifying spaces, the
twisted Morava K-homology coincides with the corresponding 
untwisted Morava K-homology. So, we see a drastic simplification 
for the following family of important spaces which arise often in the literature
(see \cite{Stolz}\cite{SSS2}\cite{SSS3}). 

\begin{corollary} [Twisted Morava K-homology of String and BString]
For the String group and its classifying space, 
we have

\vspace{-9mm}
\begin{align*}
K(5)_*\left(\op{String}; \tau_7\right) &\cong K(5)_*(\op{String}),
\\
K(6)_*\left(\op{BString}; \tfrac{1}{6}p_2\right) &\cong K(6)_*\left(\op{BString}\right),
\end{align*}
where $\tfrac{1}{6}p_2$ is 
the second fractional Pontrjagin class, classifying the fibration $K(\ZZ,7)\to \op{BFivebrane} 
\to \op{BString}$ and $\tau_7$ is its looping.
\end{corollary}

Notice that Theorem \ref{th:bu} provides us with a similar short exact sequence for connective 
covers of $\op{BU}$ and so, similarly, 
we can make the same conclusion in that case.

\medskip
\begin{theorem}
[Twisted Morava K-homology of $\op{BU}\langle k \rangle$]
Let $n$ be an odd natural number, and let 
$K(\ZZ, n+1) \to \op{BU}\langle n+2 \rangle \to \op{BU}\langle n+1 \rangle $ be a fibration 
defining a connective cover of $\op{BU}$ and $\tau_{n+2}$ the corresponding class in 
$H^{n+2}\left(\op{BU}\langle n+1\rangle; \mathbb Z \right)$. Then
\begin{align*}
K(n)_*(\op{BU}\langle n+1 \rangle; \tau_{n+2}) 
& \cong K(n)_*(\op{BU}\langle n+1 \rangle),
\\
K(n-1)_*(\op{U}\langle n+1 \rangle; \tau_{n+1}) 
& \cong K(n-1)_*(\op{U}\langle n+1 \rangle),
\end{align*}
where again the class $\tau_{n+1}$ is the looping of $\tau_{n+2}$.
\end{theorem}
%

\medskip
A more general class of spaces which will satisfy the same property is provided
by \cite[Proposition~2.0.1]{RWY}.

\medskip
Note that the restrictions on the base spaces in Theorem \ref{th:klw bo} prevent us from using the
same argument for the fibration $K(\ZZ,3)\to \op{BString} \to \op{BSpin} \to K(\ZZ, 4)$, as 
$\op{BSpin} = \underline{\op{BSpin}}_0$, and Theorem \ref{th:klw bo} requires the base to be 
$\underline{\op{BSpin}}_i$ with $i \geq 2$. However, this case can be handled using a different technique,
leading to triviality of the twisted Morava K-homology of the classifying space of the Spin group.

\begin{proposition}  
[Twisted Morava K-homology of BSpin] \label{Prop20}
We have
 \[
 K(2)_*\left(\op{BSpin}; \tfrac{1}{2}p_1\right) = 0\;,
 \]
 where $\frac{1}{2}p_1 \in H^{4}(\op{BSpin}; \ZZ)$ is the first fractional Pontrjagin class, 
 classifying the fibration $K(\ZZ,3)\to \op{BString} \to \op{BSpin}$.
\end{proposition}
\proof
By the UCT (Theorem \ref{th:uct morava}), 
we need to compute 
$K(2)_*\left(\op{BString}\right)\otimes_{K(2)_* K(\ZZ,3)} K(2)_*$. Recall also  from there that the
module structure on the second factor $K(2)_*$ is given by mapping $b_0 \in K(2)_* K(\ZZ,3)$ to $v_2$
and $b_i \in K(2)_* K(\ZZ,3)$ to $0$ for $i\geq 1$.
Now consider the exact sequence in Theorem \ref{th:klw bspin}.
From the Ravenel-Wilson computations (Theorem \ref{th:RW} above), the elements satisfy $b_i = \delta_* a_i$, 
where $K(n)_*K(\ZZ/2^j,n) \cong \bigotimes\limits_{i=0}^{j-1}R(a_i)$ and 
$R(a_k) = \ZZ/p[a_k, v_n^{\pm 1}]/(a_k^p - (-1)^{n-1}v_n^{p^k}a_k)$ for $k \geq 0$. 
So for $j=1$, we have $K(2)_*K(\Z/2, 2) \cong R(a_0) \cong \Z/2[a_0, v_2^\pm]/(a_0^2+ v_2 a_0)$, 
while $a_i=0$ for $i \geq 1$. Since $b_i=\delta_* a_i$, we have $b_0=\delta_* a_0$ while $b_i=0$ for $i \geq 1$.  
Therefore, in 
$K(2)_*\left(\op{BString}\right)\otimes_{K(2)_* K(\ZZ,3)} K(2)_*$, the relevant elements multiply as 
$$
1 \otimes v_2 = 1 \otimes b_0 = b_0 \otimes 1 = 0 \otimes 1 = 0\;.
$$
Since the element $1 \otimes v_2$ is invertible,
 the whole ring must be zero.
\endofproof

The sequence of 
 infinite loop spaces
$
K(\Z/2, 2) \to K(\Z, 3) \overset{\alpha}{\longrightarrow} \op{B String} \to \op{B Spin}
$
induces, by \cite{KLW}, the exact sequence of Hopf algebras 
\eqref{eq snake}. 
It is emphasized in \cite{KLW}\cite{KLW2} that
while, algebraically, there is a short exact sequence of Hopf algebras 
$
K(2)_*K(\Z, 3) \overset{\gamma}{\longrightarrow} 
K(2)_* \op{B String} \to K(2)_* \op{B Spin} 
$
which splits, the map $\gamma$ is {\it not} the map 
induced by $\alpha$ above. 
However, here we do not make explicit use of such splittings. 
  
\medskip
A similar argument as in the proof of Proposition \ref{Prop20} establishes Theorem 
\ref{thmvan} in the Introduction. In this case, the UCT (Theorem \ref{th:uct morava})
gives 
$
K(n)_*(B, \xi) \cong K(n)_*(E) \otimes_{K(n)_*K(\Z, n+1)}K(n)_*
$, with $K(n)_*K(\Z/2, n) \cong R(a_0) \cong \Z/2[a_0, v_n^\pm]/(a_0^2 + (-1)^n v_n a_0)$
from Theorem \ref{th:RW}. Again we have $b_0=\delta_* a_0$ while $b_i=0$ for $i \geq 1$,
and the multiplicative structure is similar.

\subsection{Twisted $K(n)$-homology of Eilenberg-MacLane spaces}
\label{Sec-EM}

We now look at bundles of Eilenberg-MacLane spaces with the base space also given by 
Eilenberg-MacLane spaces, thereby generalizing calculations of Khorami \cite{khorami}. 
We will generalize Example \ref{ExKh1} from $n=1$ to any natural number $n$. 
Note that $n$ essentially plays the role of the homotopic degree of $K(\ZZ, n+1)$ 
as well as the chromatic level of the Morava K-theory $K(n)$ being used. The proof will follow 
similar strategies to the ones taken in \cite{khorami} for the case of twisted K-homology.

\begin{theorem}
[Twisted Morava K-theory of $K(\ZZ, n)$]
Let $k: K(\ZZ, n+2) \to K(\ZZ, n+2)$ be the map induced by multiplication by 
$k$ on $\ZZ$, for $k \geq 1$. Then
\[
K(n)_*(K(\ZZ, n+2); k) = 0\;.
\]
\end{theorem}

\proof
Consider the identity map, $\op{id}: K(\mathbb Z,n+2) \to K(\mathbb Z, n+2)$, and let $P_{\rm id}$ be 
the total space of the corresponding $K(\mathbb Z, n+1)$ bundle. Notice that $P_{\rm id}$ is 
contractible by construction (it is the total space of the universal principal 
$K(\mathbb Z, n+1)$ bundle), so that $K(n)_{*}(P_{\op id}) \cong K(n)_*$.

Now consider the ``multiplication by $k$" map $\mathbb Z \xrightarrow{k} \mathbb Z$. 
Let $k:K(\mathbb Z, n+2) \to K(\mathbb Z, n+2)$ be the induced map on the Eilenberg-MacLane spaces 
and $P_k$ be the total space of the 
corresponding $K(\mathbb Z, n+1)$ bundle. Then the long exact sequence on homotopy groups for 
the principal fibration  
$K(\mathbb Z, n+1) \to P_k \to K(\mathbb Z, n+2)$ reduces to
$$
0 \longrightarrow \pi_{n+2}(P_k) \longrightarrow \pi_{n+2}\left(K(\mathbb Z, n+2)\right)\cong\mathbb Z 
\longrightarrow
\pi_{n+1}\left(K(\mathbb Z,n+1)\right)\cong \mathbb Z \longrightarrow \pi_{n+1}(P_k) \longrightarrow 0
\;,
$$
from which we observe that $P_k$ has at most two non-trivial homotopy groups. To see how the 
multiplication by $k$ fits into this picture, consider the map between $P_k$ and the universal 
$K(\mathbb Z, n+1)$ bundle
\begin{displaymath}
    \xymatrix@=1.4em{
    K(\mathbb Z, n+1) \ar[d] \ar[rr]^{\op{id}} && K(\mathbb Z,n+1) \ar[d]\\
    P_k \ar[d] \ar[rr] && \text{*} \ar[d] \\
    K(\mathbb Z, n+2) \ar[rr]^{k} && K(\mathbb Z, n+2)\;.
    }
\end{displaymath}
Here the map on base spaces is ``multiplication by $k$" by definition, and the 
map on the fibers is the identity map. This induces a map of exact sequences
\begin{displaymath}
\xymatrix{
0\ar[r] &\pi_{n+2}(P_k) \ar[d] \ar[r] & \mathbb Z \ar[d]_{k} \ar[r] &\mathbb Z \ar[d]_{\op{id}} \ar[r] 
&  \pi_{n+1}\left(P_k\right) \ar[d] \ar[r] &0\\
0\ar[r] & 0 \ar[r] & \mathbb Z  \ar[r] &\mathbb Z  \ar[r] & 0 \ar[r] & 0}
\end{displaymath}
from which we see that the top map $\mathbb Z \to \mathbb Z$ has to be multiplication by $k$. 
Consequently,  $\pi_{n+1}\left(P_k\right) = 0$ and $\pi_{n+1}\left(P_k\right) = \mathbb Z/k$, 
hence we have $P_k \simeq K(\mathbb Z/k, n+1)$.

Since $K(n)$ at prime $2$ is a 2-local theory, $K(n)_*(K(\mathbb Z/k, n+1))$ is trivial for 
$k$ odd. Together with K\"unneth isomorphism this implies that it is sufficient to look at 
$K(n)_*(K(\mathbb Z/2^j, n+1))$. From Theorem \ref{th:RW}, 
we see that 
$K(n)_*(K(\mathbb Z/2^j, n+1) \cong K(n)_*$
and, therefore, for all $n\geq 1$ we have
\[
K(n)_*\big(K(\mathbb Z, n+2); k\big) \cong K(n)_* \tensor_{K(n)_{*}(K(\mathbb Z, n+1)} K(n)_*\;.
\]
Now recall that, in the module structure on the second factor, $b_0$ from $K(n)_*K(\mathbb Z, n+1)$ 
is mapped to 1. On the other hand, $b_0$ maps to 0 in $K(n)_*(K(\mathbb Z/2^j, n+1) \cong K(n)_* $,
so that we finally arrive at $K(n)_*(K(\mathbb Z, n+2); k) = 0$, for any $k>0$.
\endofproof

\vspace{-5mm} 
\subsection{Twisted $K(n)$-(co)homology of spheres}
\label{Sec-Sph} 

We now generalize  Example \ref{ExKh2} 
from $n=1$ to any natural number $n$, which plays the role of the dimension of the sphere (minus 2) 
as well as the chromatic level. This time 
the results will differ and the proof will depart drastically from 
those of \cite{khorami};  instead will use the twisted Atiyah-Hirzebruch 
spectral sequence (AHSS) of \cite{SW}, 
i.e., Theorem \ref{th:ahss}. We find that the twisted Morava K-theory is 
essentially trivial, with the meaning of triviality depending on chromatic level. 
For levels greater than 1, it is given by the underlying integral cohomology tensored 
with the coefficients of the theory, while for level 1, it is trivial in the 
sense of being actually zero. 

\begin{theorem}
[Twisted Morava K-theory of spheres]
Fix an integer $n>1$ and let $\sigma_{n+2}$ be the generator of the group 
$H^{n+2}(S^{n+2}; \ZZ) = \left[S^{n+2}, K(\ZZ, n+1)\right]\cong \ZZ$. Then the twisted 
$n$th Morava 
K-cohomology of the $(n+2)$-sphere with  twist $\tau_{n+2}$ 
given by any multiple of $\sigma_{n+2}$ is given 
by cohomology tensored with the coefficient ring:
 \[
 K(n)^{*}\left(S^{n+2}; \tau_{n+2}\right) = H^*(S^{n+2}; \Z) \otimes K(n)_*\;.
 \]
\end{theorem}
\proof
Let $P$ again denote the total space of the corresponding $K(\ZZ, n+1)$ bundle. Recall from 
Theorem \ref{th:ahss} that the first non-trivial differential in the twisted AHSS for 
Morava K-theory is 
\[
d_{2^n-1}(x) = Q_n(x) + Q_{n-1}\cdots Q_1(H) \cup x\;.
\]
Because the differential is a module homomorphism, is suffices to consider cases $x = 1$ 
and $x$ dual to the fundamental class of $S^{n+2}$. In both cases $Q_n(x) = 0$ by 
dimension, since $H^{*}(S^{n+2}, K(n)^{*})$ concentrated in degree $n+2$, so the 
target of $Q_{n}$ is trivial. 

When $n>1$, we do have the part $Q_{n-1}\cdots Q_1(H)$ in the second term. Then,
 by the same dimension argument, $Q_{n-1}\cdots Q_1(H) = 0$ since the target of 
 $Q_{n-1}\cdots Q_1$ is of higher degree cohomology. Therefore, the first possibly 
 non-trivial differential is actually trivial. Also note that all the subsequent differentials 
 must vanish too (they are even longer so will land in even higher cohomology groups). 
 Since this holds for every $x$, the spectral sequence collapses. 
\endofproof

The case $n=1$ can be handled separately, by either using the non-twisted AHSS to 
compute $K(n)_*P$, or using the computations of twisted K-theory 
\cite{khorami} together with Theorem \ref{th:mod p} : by the computations in Section 
\ref{sec:khorami}, twisted K-homology for $S^3$ twisted by the generator 
$\sigma_3\in H^{3}(S^3;\ZZ)$ is $K^{\sigma_{3}}_*(S^3) = 0$. But, from 
Theorem \ref{th:mod p}, twisted $K(1)$ fits into the short exact sequence
\vspace{-2mm} 
\[
0 \longrightarrow K_n^{\sigma_3}(S^3) \otimes \ZZ/2 \longrightarrow 
 K(1)_n(S^3;\sigma_3) \longrightarrow  
 \op{Tor}_{1}\big(K_{n-1}^{\sigma_3}(X), \ZZ/2\big)\;. 
 \]

 \vspace{-2mm} 
\noindent Since $K^{\sigma_{3}}_*(S^3) = 0$, both the first and the third term of this exact sequence 
are zero, so that the middle term is zero as well. On the other hand, starting with a twist 
which is a multiple of the generator we can use $K_*(S^3; n \sigma_3)=\Z/n\Z$. 
Therefore, using the relation between Morava $K(1)$-homology at the prime 2 and 
K-homology, we arrive at the following. 

\begin{proposition}
[Twisted first Morava $K$-homology of the 3-sphere]
\label{thm 3-sph}
The Morava K(1)-homology at $p=2$ (i.e. mod 2 K-homology) 
of the 3-sphere with a twist a multiple of the generator $\sigma_3$ is  
 $$
    K(1)_{*}(S^{3}; n\sigma_3) = \Z/(2, n)[v_1^{\pm 1}]\;,
$$
while it vanishes for $n=1$,
$
    K(1)_{*}(S^{3}; \sigma_3) = 0
$.
\end{proposition}

This is the mod 2 version of Khorami's theorem \cite{khorami}  -- see Example \ref{ExKh2}.
Notice that this is consistent with the pattern we observed before: for `nice enough' spaces, 
twisted Morava K-theory groups are either zero or equal to untwisted groups.

%


\subsection{Twists by mod 2 Eilenberg-MacLane spaces and torsion connective covers}
\label{Sec-mod2}

We would like to complete our investigations of connective covers of $\op{BO}$. 
So far we have been focusing solely on those covers which can be viewed as bundles 
of integral Eilenberg-MacLane spaces, i.e., those levels of the Whitehead tower of 
$\op{BO}$ in diagram \eqref{eq:tower} which have maps to $K(\ZZ, m)$. We would like 
to perform a similar analysis for the remaining ``non-integral" covers, for instance,
the analogue of orientation 
$K(\ZZ/2, 8) \to \op{BO}\langle 10\rangle = \op{B2\textit{-}Orient} \to 
\op{BO}\langle 9\rangle = \op{BFivebrane}$
and the analogue of Spin structure 
$K(\ZZ/2, 9) \to\op{BO}\langle 11\rangle = \op{B2\textit{-}Spin} 
\to \op{BO}\langle 10\rangle = \op{B2\textit{-}Orient}$, 
which are $K(\ZZ/2, m)$-bundles. However, so far we have been lacking the definition of 
Morava K-theory twisted by non-integral Eilenberg-MacLane spaces. The purpose of this 
section is to fill that gap. 

\medskip
Instead of focusing solely on $p=2$, we will discuss twists of $K(n)$ by 
$K(\ZZ/p^j, m)$ for all primes $p$ and $j \geq 1$. 
From the description of Morava K-theory in Theorem \ref{th:RW} we see $K(n)_{*}K(\ZZ/p^j,n)$ is one of the factors 
of $K(n)_{*}K(\ZZ, n+1)$, and $K(n)_{*}K(\ZZ/p^j,m) = K(n)_{*}K(\ZZ,m) = K(n)_*$ for $m>n$. Therefore, we should 
expect a similar theory as for twists by integral Eilenberg-MacLane spaces. In fact, the proofs in \cite{SW} 
transport to the mod $p$ case with little to no modification, and so we only outline them.

\medskip
Recall from \cite{ABG}\cite{ABGHR} that a twist of theory $R$ by a space $Y$ is an element of 
$[Y, \op{BGL}_1 R]$. The following fact provides us with an obstruction-theoretic way to classify 
these maps.

\begin{proposition}[{\cite[Proposition~1.6]{SW}}]\label{th:obstruction}
Let $R$ be an $A_{\infty}$ ring spectrum, $Z = \Omega X$ and $R_*(Z)$ is flat over $R_*$. If the 
obstruction groups 
   $ \op{Ext}{}_{R_*(Z)^{\rm op}}^k(R_*, \Omega^s R_*)$
vanish for $s=k-1,k-2$ and any $k \geq 1$, then there is a bijection
$
\op{tw}_R(X) \leftrightarrow \op{Hom}_{\op{R_*{\textbf{--}alg}}}(R_*(Z), R_*)
$.
Moreover, the obstruction groups lie in the $E_2$-term of the cobar spectral sequence
\begin{equation}
\label{eq:ss}
\op{Ext}{}_{R(Z)^{\rm op}}^k(R_*, \Omega^s R_*) \; \Longrightarrow \; E^{k-s}(BZ)\;.
\end{equation}
\end{proposition}

\medskip
Notice that when $R = K(n)$, the flatness requirement is automatically satisfied for 
any $Z$, since any $K(n)_*(Z)$ is free over $K(n)_*$. 
Now we can establish the following mod 2 analogue of Theorem \ref{SW-tw}.

\begin{theorem}
[Morava K-theory twisted by mod $p$ Eilenberg-MacLane spaces] We have:
	 \vspace{-2mm}
	    \item {\bf (i)} There are no non-trivial twists of $K(n)$ by $K(\ZZ/p^{j}, m)$  for $m >n+1$;
	    \vspace{-2mm} 
	    \item {\bf (ii)} There are no non-trivial twists of $K(n)$ by $K(\ZZ/p^j, n+1)$ at $p \neq 2$;
	    \vspace{-2mm} 
		\item {\bf (iii)} $\op{tw}_{K(n)}(K\left(\ZZ/2^j, n+1\right) \cong \ZZ/2^j$. 
	\end{theorem}

\begin{remark}
Notice the shift in degree compared to integral Eilenberg-MacLane spaces. It is the 
same shift in degree that occurs in Theorem \ref{th:RW}.
\end{remark}
\proof (Outline)
We will use Proposition \ref{th:obstruction} with $X = K(\ZZ/p^j, m)$, $Z = K(\ZZ/p^j, m-1)$, and $R=K(n)$. 
From Theorem \ref{th:RW}, if $m>n+1$ then $K(n)_*(K(\ZZ/p^j, m-1)) = K(n)_*$. Consequently,  the 
obstruction group is $\op{Ext}_{R_*}^k(R_*, \Omega^s R_*) = 0$, so that the twists are given as 
$$
\op{tw}_{K(n)}(K(\ZZ/2, m)) = \op{Hom}_{K(n)_*{-{\rm alg}}}(K(n)_*, K(n)_*) = \{\ast\} .
$$
Just as in the integral case, the spectral sequence  \eqref{eq:ss} collapses by the work 
of Ravenel and Wilson \cite{RW}, and the obstruction groups vanish, leading to
$$
\op{tw}_{K(n)}(K(\ZZ/p^j, n+1)) = \op{Hom}_{K(n)_*{\rm -alg}}
\big(K(n)_*K(\ZZ/p^j, n+1), K(n)_*\big) .
$$ 
Since twistings of $K(\Z, n+2)$ restrict (along the Bockstein) to twistings of $K(\Z/p^j, n+1)$, 
this is a quotient of  
$\op{Hom}_{K(n)_*{\rm -alg}}\big(K(n)_*K(\ZZ,n+2), K(n)_*\big) = \op{tw}_{K(n)}(K(\ZZ, n+2))$,
the twisting space for $K(\Z, 2)$ from Theorem \ref{th:RW}. But we know that, for $p>2$,
 the latter is trivial.

Now fix $p=2$, and recall from Theorem \ref{th:RW} 
that $K(n)_*K(\ZZ/p^j,n) \cong \bigotimes\limits_{i=0}^{j-1}R(a_i)$. As in the proof of \cite[Theorem~3.3]{SW}, 
$\op{Hom}_{K(n)_*{\rm -alg}}\big(K(n)_*(K(\ZZ/2^j, n+1)), K(n)_*\big)$ 
is determined by the images of the elements $a_i$, for $0\leq i \leq j-1$. By degree reasons, there is only 
one possible target for each $a_i$ in the coefficient ring 
$K(n)_*$. So an element of 
$\op{Hom}_{K(n)_*{\rm -alg}}(K(n)_*(K(\ZZ/2^j, n+1)), K(n)_*)$ 
is determined by the $j$ elements among 
the $a_i$ which are mapped to zero, and there are $2^j$ elements. By identifying  
$\op{Hom}_{K(n)_*{\rm -alg}}\big(K(n)_*(K(\ZZ/2^j, n+1)), K(n)_*\big)$ with a 
subset of $K(n)^{*}[x]/x^{2^j}$ (closed under multiplication, being the group-like
elements of the Hopf algebra)  it 
is possible to obtain a group structure on it.
\endofproof

This allows us to seek direct analogues of the constructions in \cite{SW}, as recalled in Section \ref{Sec-tmo}. 
In particular, since $\op{tw}_{K(n)}(K\left(\ZZ/2, n+1\right) \cong \ZZ/2$, 
we can present an analogue of Definition \ref{def-u}. 

\begin{definition}
{\bf (i)}    The \textit{universal twist} of $K(n)$ by the mod 2 Eilenberg-MacLane space 
$K(\ZZ/2,n+1)$ is the non-zero element of $\op{tw}_{K(n)}(K\left(\ZZ/2, n+1)\right)$.
\vspace{-2mm} 		
\item {\bf (ii)} Let $h \in H^{n}(X; \ZZ/2)$ be a mod 2 cohomology class. Then 
\emph{Morava K-theory of $X$ twisted by $h$} is defined to be  
$K(n)_{*}(X; h):= K(n)^{u(h)}_*(X)$.
\end{definition}
The universal coefficient theorem analogue of Theorem \ref{th:uct morava} is also true in this 
case, and the proof follows the proof of that theorem with obvious changes (hence we omit it to avoid repetition).

\begin{theorem}
[UCT for Morava K-theory twisted by mod $2$ EM spaces]
If $h \in H^{n}(X; \ZZ/2)$, and $P$ denotes the total space of the bundle classified by $h$, 
then 
$$
K(n)_{*}(X; h) \cong K(n)_{*}\left(P\right) \otimes_{K(n)_{*}K\left(\ZZ/2, n\right)} K(n)_{*}\;. 
$$
\end{theorem}

\medskip
Equipped with this result, we can conclude our investigation of the Whitehead tower of $\op{BO}$.
In diagram \eqref{eq:tower}, the notation introduced in \cite{9} for the $\ZZ/2$'s in the Whitehead tower 
refers to the Bott periodicity cycles of length 8. That is, 
${\rm B}(m\textit{-}{\rm Orient})$ and ${\rm B}(m\textit{-}{\rm Spin})$ are those mod 8 analogues of orientation 
and Spin structure that occur in Bott periodicity cycle $m$. So we have a sequence 
${\rm B}(1\textit{-}{\rm Orient})={\op{B{\rm SO}}}$ and $\op{B(1\textit{-}Spin)}=\op{BSpin}$ in the first
cycle, $\op{B(2\textit{-}Orient}):= {\rm BO} \langle 10 \rangle$ and $\op{B(2\textit{-}Spin)}:={\rm BO} \langle 11 \rangle$
in the second cycle and so on (see the Whitehead tower, Diagram \eqref{eq:tower}). 

\medskip
Applying Theorem \ref{th:klw bo} for $G = \op{BSO}_i$ and $G=\op{BO}_i$ (with $i \geq 2)$, 
we get the following mod 2 version of Theorem \ref{Kn-BO}, with a basically identical proof. 

\begin{theorem}
[Twisted Morava K-theory for  ${\rm B}(m\textit{-}{\rm Orient})$ and ${\rm B}(m\textit{-}{\rm Spin})$ structures]
Let $\op{BO}\langle n \rangle$ be a connective cover of $\op{BO}$ with 
$n  = 1 \mod \,8$ or $n= 2 \mod \,8$, and let $h_{n+1}$ be the class in 
$H^{n+1}(\op{BO}\langle n \rangle ;\mathbb Z/2)$ classifying the connective cover fibration. 
Then
\begin{align*}
	K(n)_*\left( \op{BO}\langle n \rangle; h_{n+1} \right) &\cong
	K(n)_*\left( \op{BO}\langle n \rangle\right) \quad \text{ for any }n \geq 8\;,
	\\
	K(n)_*\left( \op{O}\langle n \rangle; h_{n} \right) &\cong
	 K(n)_*\left( \op{O}\langle n \rangle\right) \quad \text{ for any }n \geq 7\;,
\end{align*}
where $h_{n} \in H^n(O \langle n \rangle; \Z/2)$ is the looping of $h_{n+1}$. 
\end{theorem}

The cohomology groups $H^{n+1}(\op{BO}\langle n \rangle ;\mathbb Z/2)$, in which our
twists live, have been determined by Stong \cite{Stong}. They are 
isomorphic to $H^{n+1}\big(K(\pi_n(\op{BO}), n); \Z/2\big)/I \otimes \Z/2 [\theta_i | L(i) > \phi(n)]$, where
$I$ is some ideal depending on $n$, $\theta_i$ is a class in $H^i(\op{BO}; \Z/2)$, and 
$L(i)$ and $\phi(n)$ are integers constructed from $i$ and $n$, respectively. 
Specific values can be read off from the above paper of Stong. 

\medskip
Note that the only connective covers of $\op{O}$ and $\op{BO}$ that we have not investigated so far are  
$\op{Spin}$, $\op{SO}$, and $\op{BSO}$. The first two are not directly interesting for our purposes: 
$\op{Spin}$ is defined via a map to $K(\ZZ/2,1)$ and $\op{SO}$ is  defined via a map to $K(\ZZ/2, 0)$. 
This would mean that the corresponding twisted Morava K-theory  has to be at height 0, i.e., rational 
cohomology, for which the result is classical.  
This is the statement that $H^*(\op{SO}; t)\cong H^*(\op{Spin}; \Z)$, where
$t \in H^1(\op{SO}; \Z/2)$ classifies the double cover $\op{Spin} \to \op{SO}$.
For $B{\rm SO}$, the Whitehead tower of the orthogonal group, Diagram \eqref{eq:tower},
gives us the fibration
\[
K(\ZZ/2, 1) \longrightarrow \op{BSpin} \longrightarrow \op{BSO} \xrightarrow{w_2} 
K(\ZZ/2, 2)\;,
\]
where $w_2$ is the second Stiefel-Whitney class. However, as shown in \cite[Section~5.3]{KLW}, 
the induced map $K(n)_{*} (K(\ZZ/2,1)) \to K(n)_{*} (\op{BSpin})$ has to be trivial, 
so it sends $b_0$ to 0 in $K(n)_*(\op{BSpin})$. Therefore,  by the same argument as in Theorem 
\ref{th:klw bspin}, we have the following at chromatic level one.

\begin{corollary}
[Mod 2 twisted Morava $K(1)$ of BSO at the prime 2]
With the twist given by the second Stiefel-Whitney class, we have
$$
K(1)_*(\op{BSO}; w_2) \cong 0\;.
$$
\end{corollary} 

Throughout, we have only considered specific twists. It would be interesting to
consider other types of twists or (further) arbitrary multiples of the twisting class. It would 
also be interesting to consider groups and their connective covers unstably. These will 
undoubtedly be subtle, as witnessed by the computations in \cite{Rao1}\cite{Rao}\cite{Ni} 
of Morava K-theory of $\op{SO}(m)$ and $\op{Spin}(m)$ and by \cite{SSh} in studying  
connective covers of these unstable groups.


\end{document}